\documentclass{article}
\usepackage{amssymb,amsmath,bm,geometry,graphics,graphicx,url,color}

\def\no{\noindent}
\def\pmatrix{\left(\begin{array}}
\def\endpmatrix{\end{array}\right)}

\newtheorem{theo}{Theorem}

\newtheorem{defi}{Definition}

\newtheorem{assum} {Assumption}

\newtheorem{prop} {Proposition}

\newtheorem{rem}{Remark}

\newcommand{\norm}[1]{\left\Vert#1\right\Vert}

\title{Long-time  energy analysis of  extended RKN
integrators for  muti-frequency highly oscillatory Hamiltonian
systems}

\author{Bin Wang\,
\footnote{School of Mathematical Sciences, Qufu Normal University,
Qufu 273165, P.R. China; Mathematisches Institut, University of
T\"{u}bingen, Auf der Morgenstelle 10, 72076 T\"{u}bingen, Germany.
The research is supported in part by the Alexander von Humboldt
Foundation and by the Natural Science Foundation of Shandong
Province (Outstanding Youth Foundation) under Grant ZR2017JL003.
E-mail:~{\tt wang@na.uni-tuebingen.de} } \and Xinyuan
Wu\thanks{School of Mathematical Sciences, Qufu Normal University,
Qufu 273165, P.R. China; Department of Mathematics, Nanjing
University, Nanjing 210093, P.R. China. The research is supported in
part by the National Natural Science Foundation of China under Grant
11671200. E-mail:~{\tt xywu@nju.edu.cn}} }

\begin{document}
\maketitle

\begin{abstract} In this paper, we study the long-time  near conservation of the
total and oscillatory energies for extended RKN  (ERKN) integrators
when applied to muti-frequency highly oscillatory Hamiltonian
systems. We consider one-stage explicit symmetric integrators and
show their long-time behaviour of numerical energy conservations by
using modulated multi-frequency Fourier expansions. Numerical
experiments are carried out  and the numerical results demonstrate
the remarkable long-time  near conservation of the energies for the
ERKN integrators and support our theoretical analysis presented in
this paper.
\medskip

\no{\bf Keywords:}  Long-time  energy conservation \and Modulated
Fourier
 expansions \and
Muti-frequencies highly oscillatory systems \and  Hamiltonian
systems \and Extended RKN integrators

\medskip
\no{\bf MSC:}65P10 \and 65L05

\end{abstract}

\section{Introduction}
The study of numerical energy preservation is an important aspect of
numerical analysis in the sense of structure-preserving algorithms
when applied to  Hamiltonian systems. This paper is devoted to
muti-frequency highly oscillatory Hamiltonian systems with the
following Hamiltonian function
 \begin{equation}\label{H}
\begin{array}
[c]{ll}%
H(q,p) &
=\dfrac{1}{2}\textstyle\sum\limits_{j=0}^{l}\Big(\norm{p_j}^{2}+\dfrac
{\lambda_j^{2}}{\epsilon^2}\norm{q_j}^{2}\Big)+U(q),
\end{array}
\end{equation}
where $q=(q_{0},q_1,\ldots,q_l),\ p=(p_{0},p_1,\ldots,p_l)$ with
$q_j,\ p_j\in \mathbb{R}^{d_j}$, $\lambda_0=0$ and $\lambda_j\geq1$
for $j\geq1$ are distinct real numbers , $\epsilon$ is a small
positive parameter, and $U(q)$ is a smooth potential function. We
pay attention to the muti-frequency case where $l>1.$  As is known,
this system has the oscillatory energy of the $j$th frequency as
\begin{equation*}I_j(q,p)=\dfrac{1}{2} \Big(\norm{p_j}^{2}+\dfrac
{\lambda_j^{2}}{\epsilon^2}\norm{q_j}^{2}\Big),
\end{equation*}
 and its total oscillatory energy is given by
$
I (q,p)=\sum\limits_{j=1}^{l}I_j(q,p). 
$  Letting
\begin{equation*}
\begin{aligned}&\lambda=(\lambda_1,\ldots,\lambda_l),\quad
k=(k_1,\ldots,k_l),\quad k\cdot
\lambda=k_1\lambda_1+\cdots+k_l\lambda_l,
\end{aligned}
\end{equation*}
and denoting the resonance module by
\begin{equation}\mathcal{M}= \{k\in \mathbb{Z}^{l}:\ k\cdot
\lambda=0\},
\label{M}%
\end{equation}
it follows from  the analysis in \cite{Benettin89} that the
quantities
\begin{equation*}
\begin{aligned}
& I_{\mu}(q,p)=  \sum\limits_{j=1}^l \frac{\mu_j}{\lambda_j}
I_j(q,p)\quad \mathrm{with}\ \mu\ \mathrm{orthogonal\ to}\
\mathcal{M}
\end{aligned}
\end{equation*}
are approximately preserved under a diophantine non-resonance
condition outside $\mathcal{M}$. Here,  it is clear that
$I_{\lambda}=I(q,p)$.

This kind of  muti-frequency  highly oscillatory systems often
arises in various fields such as applied mathematics, molecular
biology, astronomy, and classical mechanics (see, e.g.
\cite{hairer2006,Hochbruck2010,Wang2018-1,wu2013-book}). In  recent
years, many effective numerical methods have been developed and see,
e.g.
\cite{Grimm2006,Hairer09,Hochbruck1999,Li_Wu(sci2016),wang-2016,wang2017-ANM,wu2017-JCAM}
as well as the references contained therein. In \cite{wu2010-1}, the
authors  formulated  a kind of trigonometric integrators called as
  extended Runge--Kutta--Nystr\"{o}m (ERKN) integrators for solving
muti-frequency  highly oscillatory systems. Some important
properties of these integrators were further studied  in
\cite{wu2013-ANM,BIT_Wu,wu2013-book}. Very recently, the long-time
 energy conservation of ERKN integrators for highly
oscillatory Hamiltonian systems with one frequency was researched in
\cite{17-new}.   On the basis of this work, this paper is devoted to
the numerical energy analysis of ERKN integrators for muti-frequency
highly oscillatory Hamiltonian systems.

For the analysis of energy preservation,  modulated Fourier
expansions are an elementary and useful analytical tool. It was
firstly developed in \cite{Hairer00} and then was used as an
important mathematical tool in studying the long-time behaviour of
numerical methods for differential equations (see, e.g.
\cite{Cohen06,Cohen12,Cohen15,Cohen03,Cohen05,Cohen06-1,Gauckler10-1,Hairer12-1,Hairer16,iserles08,McLachlan14,Sanz-Serna09,Stern09,17-new}).
The long-time analysis of some trigonometric integrators for
muti-frequency oscillatory Hamiltonian systems has been given in
\cite{Cohen05}. In this paper we extend the long-time energy
preservation results of \cite{Cohen05,17-new} to the ERKN
integrators for multi-frequency cases. As is known, resonance
frequencies may exist for multi-frequency oscillatory Hamiltonian
systems. Hence, compared with the analysis of one-frequency case in
\cite{17-new}, a new and important aspect of muti-frequency case is
possible resonance among the frequencies, which is   similar to the
analysis made in \cite{Cohen05}.

The remainder of this paper is organised as follows. In Section
\ref{sec preliminaries}, we briefly summarise   ERKN integrators for
the muti-frequency   Hamiltonian systems \eqref{H} and present some
preliminaries.
 The modulated Fourier expansion of ERKN integrators are
derived and analysed in Section \ref{sec:Analysis of the methods}
and   two almost-invariants of the modulated Fourier expansions  are
studied in Section \ref{sec:Almost-invariants}. Then Section
\ref{sec:Long-time near-conservation} presents the main result
concerning the long-time near energy conservation. Numerical
experiments are accompanied  in Section \ref{sec:examples}. The last
section is concerned with the conclusions of this paper.
\section{Preliminaries}\label{sec preliminaries}
\subsection{ERKN integrators} \label{subsec:Formulation}

Rewrite the highly oscillatory system \eqref{H}    as a system of
second-order differential equations
\begin{equation}
q''=-\Omega^2 q+g(q), \qquad
q(0)=q^0,\ \ p(0)=p^0, \label{prob}%
\end{equation}
where  $\Omega=\textmd{diag}(\omega_0
I_{d_0},\omega_1I_{d_1},\ldots, \omega_lI_{d_l})$ with
$\omega_j=\lambda_j/\epsilon$ and $g(q)=-\nabla U(q).$ A kind of
trigonometric integrators called as ERKN integrators has been
  developed (see, e.g. \cite{wu2010-1}), and  the one-stage ERKN explicit scheme  will be discussed in detail in this paper.

\begin{defi}
\label{erkn}  (See \cite{wu2010-1}) A one-stage  explicit ERKN
integrator for \eqref{prob} is defined by%
 \begin{equation}
\begin{array}
[c]{ll}%
Q^{n+c_{1}} &
=\phi_{0}(c_{1}^{2}V)q^{n}+hc_{1}\phi_{1}(c_{1}^{2}V)p^{n},\\
q^{n+1} & =\phi_{0}(V)q^{n}+h\phi_{1}(V)p^{n}+h^{2}
 \bar{b}_{1}(V)g(Q^{n+c_{1}}),\\
p^{n+1} & =-h\Omega^2\phi_{1}(V)q^{n}+\phi_{0}(V)p^{n}+h\textstyle
b_{1}(V)g(Q^{n+c_{1}}),
\end{array}
  \label{methods}%
\end{equation}
where   $h$ is a stepsize, $c_1$ is real constant satisfying $0\leq
c_1\leq1$ , $b_{1}(V)$ and $\bar{b}_{1}(V)$  are matrix-valued and
uniformly bounded functions of $V\equiv h^{2}\Omega^2$, and
\begin{equation}
\phi_{j}(V):=\sum\limits_{k=0}^{\infty}\dfrac{(-1)^{k}V^{k}}{(2k+j)!},\qquad j=0,1,\ldots.%
\label{Phi01}%
\end{equation}
\end{defi}

The following three results of ERKN integrators will be useful in
this paper.

\begin{theo}
\label{order conditions} (See \cite{wu2010-1,wu2013-book}) The
one-stage explicit ERKN integrator \eqref{methods}  is of order two
if and only if
\[%
\begin{array}
[c]{l}%
 b_1(V)=\phi_{1}(V)+\mathcal{O}(h^2),\ \
 c_1b_1(V)=\phi_{2}(V)+\mathcal{O}(h),\ \
\bar{b}_1(V)=\phi_{2}(V)+\mathcal{O}(h).
\end{array}
\]
\end{theo}

\begin{theo}\label{symmetric thm}(See \cite{wu2013-book})
If and only if
\begin{equation}\begin{aligned}\label{sym cond}&c_1=1/2,\
\bar{b}_{1}(V)=\phi_{1}(V)b_{1}(V)-\phi_{0}(V)\bar{b}_{1}(V),\
\phi_{0}(c_{1}^{2}V)\bar{b}_{1}(V)=c_{1}\phi_{1}(c_{1}^{2}V)b_{1}(V),
\end{aligned}\end{equation}
then the one-stage explicit ERKN integrator \eqref{methods} is
symmetric.
\end{theo}

\begin{theo}\label{symplectic thm} (See \cite{wu2013-book})
If there exists a real number $d_1$ such that
\begin{equation*}\begin{aligned}
&\phi_{0}(V)b_{1}(V)+V\phi_{1}(V)\bar{b}_{1}(V)=d_{1}\phi_{0}(c_{1}^{2}V),\ \ \ d_1\in \mathbb{R},\\
&\phi_{1}(V)b_{1}(V)-\phi_{0}(V)\bar{b}_{1}(V)=c_{1}d_{1}\phi_{1}(c_{1}^{2}V),
\end{aligned}\end{equation*}
then the one-stage explicit ERKN integrator \eqref{methods} is
symplectic.
\end{theo}

\subsection{Notations} \label{subsec:Notations}
Let $V=h^{2}\Omega^2$. It follows from \eqref{Phi01} that
\begin{equation*}
\phi_{0}(V)=\cos(h\Omega),\qquad  \phi_{1}(V)=\textmd{sinc}(h\Omega):=(h\Omega)^{-1}\sin(h\Omega).
\end{equation*}
Throughout  this paper, we use the notations $\bar{b}_{1}(h\Omega)$
and $ b_{1}(h\Omega)$ to denote the coefficients appearing in the
ERKN method \eqref{methods}. Moreover, we  also adopt the following
notations which  appeared in \cite{Cohen05}:
\begin{equation*}
\begin{aligned}
&\omega=(\omega_1,\ldots,\omega_l),\quad \langle
j\rangle=(0,\ldots,1,\ldots,0),\quad |k|=|k_1|+\cdots+|k_l|.
\end{aligned}
\end{equation*}
For  the resonance   module \eqref{M}, we let $\mathcal{K}$ be a set
of representatives of the equivalence classes in
$\mathbb{Z}^l\backslash \mathcal{M}$ which are chosen such that for
each $k\in\mathcal{K}$ the sum $|k|$ is minimal in the equivalence
class $[k] = k +\mathcal{M},$ and that with $k\in\mathcal{K}$, also
$-k\in\mathcal{K}.$ We denote, for the positive integer $N$,
\begin{equation}\mathcal{N}=\{k\in\mathcal{K}:\ |k|\leq N\},\ \ \ \ \ \mathcal{N}^*=\mathcal{N}\backslash
\{(0,\ldots,0)\}.
\label{mathcalN}%
\end{equation}
 In this paper, we  use the following  operator  which has
been defined in \cite{hairer2006}
\begin{equation*}
\begin{aligned}\mathcal{L}(hD):&=\mathrm{e}^{hD}-2\cos(h\Omega)+\mathrm{e}^{-hD}=2\big(\cos(\mathrm{i}
hD)-\cos(h\Omega)\big)\\
&= 4\sin\big(\frac{1}{2}h\Omega+\frac{1}{2}\mathrm{i}
hD\big)\sin\big(\frac{1}{2}h\Omega-\frac{1}{2}\mathrm{i}
hD\big),\end{aligned}
\end{equation*}
where $D$ is the differential  operator. It is easy to verify that
$(hD)^m x(t)=h^m x^{(m)}(t)$ for $m=0,1,\ldots,$ and
$\mathrm{e}^{hD} x(t)=  x (t+h).$

We consider the application of such an operator to functions of the
form $\mathrm{e}^{\mathrm{i}(k \cdot \omega) t}.$ By Leibniz' rule
of calculus, one has
$$(hD)^m \mathrm{e}^{\mathrm{i}(k \cdot \omega) t} z(t)=
\mathrm{e}^{\mathrm{i}(k \cdot \omega) t}(hD+\mathrm{i}(k \cdot
\omega) h)^m z (t),$$  which yields $ f(hD)
\mathrm{e}^{\mathrm{i}(k \cdot \omega) t} z(t)=
\mathrm{e}^{\mathrm{i}(k \cdot \omega) t}f(hD+\mathrm{i}(k \cdot
\omega) h)  z (t),
$ where $$f(hD+\mathrm{i}(k \cdot \omega) h)z
(t)=\sum\limits_{m=0}^{\infty}\frac{f^{(m)}( \mathrm{i}(k \cdot
\omega) h)}{m!}h^m z^{(m)}(t).$$ Furthermore, we have the following
proposition  of the operator.

\begin{prop}\label{lhd pro}
The   Taylor expansions of $\mathcal{L}(hD)$ and
$\mathcal{L}(hD+\mathrm{i}(k \cdot \omega) t)$ are
\begin{equation*}
\begin{aligned}\mathcal{L}(hD)=&4\sin^2(h \Omega/2)-I(\mathrm{i} h
D)^2+\ldots,\\
\mathcal{L}(hD+\mathrm{i}(k \cdot \omega) h)=&\big(2\cos((k \cdot
\omega) h)I-2\cos(h \Omega)\big)+2\sin((k \cdot \omega)
h)I(\mathrm{i}
h D)\\
&-\cos((k \cdot \omega) h)I(\mathrm{i} h D)^2+\ldots.
\end{aligned}
\end{equation*}
\end{prop}

\section{Modulated Fourier expansion of the integrators} \label{sec:Analysis of the methods}
Before presenting the analysis of  long-time conservation, we make
the following assumptions.  The first four assumptions have been
considered in \cite{Cohen05}.
\begin{assum}\label{ass}
 $\bullet$   The initial values are assumed to satisfy
\begin{equation*}
\frac{1}{2} \norm{p^0}^2+\frac{1}{2} \norm{\Omega q^0}^2\leq E. 
\end{equation*}

 $\bullet$ It is assumed that the numerical solution $Q^{n+c_{1}}$  stays in a compact set on which the potential $U$ is
 smooth.

 $\bullet$ A lower bound  is posed  for  the stepsize
$
h/\epsilon \geq c_0 > 0. 
$

 $\bullet$ Assume that the following numerical non-resonance condition  holds
\begin{equation}
|\sin(\frac{h}{2 \epsilon}(k\cdot \lambda))| \geq c \sqrt{h}\ \
\mathrm{for} \ \ k \in \mathbb{Z}^l\backslash \mathcal{M}
   \ \   \mathrm{with} \ \  |k|\leq N\label{numerical non-resonance cond}%
\end{equation}
for some $N\geq2$ and $c>0$. In this paper, the $\mathcal{N}$
 given in \eqref{mathcalN} is defined for this   $N$.

 $\bullet$ The ERKN integrators are required to satisfy the symmetry conditions \eqref{sym cond}.
 Moreover, it is assumed that  $$ |b_1(h\omega_j)|\leq C_2
 |\textmd{sinc}(h\omega_j/2)|,$$
for $j=1,\ldots,l.$
\end{assum}

 \begin{rem} It is clear that we consider the numerical
non-resonance condition \eqref{numerical non-resonance cond} in the
analysis of this paper, which is the same as that in \cite{Cohen05}.
We also noted that the long-term analysis of some integrators for
oscillatory  systems under minimal non-resonance conditions has
recently been presented in \cite{Cohen15}. The long-time analysis of
ERKN integrators under minimal non-resonance conditions will be our
next work in the near future. \end{rem}

We will establish   a modulated Fourier expansion for the ERKN
integrators by the following theorem. It is the multi-frequency
version of \cite{17-new}. Its proof follows the lines of the proof
of the corresponding theorem given in \cite{17-new} but with rather
obvious adaptations. In the proof of this theorem,  we just briefly
highlight the main differences and ignore  the same derivations for
brevity.

\begin{theo}\label{energy thm}
Suppose that  Assumption \ref{ass} is true.     The ERKN integrator
\eqref{methods} admits the expansions
\begin{equation*}
\begin{aligned} &q^{n}= \zeta(t)+\sum\limits_{k\in\mathcal{N}^*} \mathrm{e}^{\mathrm{i}(k \cdot \omega) t}\zeta^k(t)+R_{h,N}(t),\\
&p^{n}= \eta(t)+\sum\limits_{k\in\mathcal{N}^*} \mathrm{e}^{\mathrm{i}(k \cdot \omega) t}\eta^k(t)+S_{h,N}(t),\\
\end{aligned}
\end{equation*}
for $0 \leq t=nh \leq T$. The remainder terms are bounded by
\begin{equation}
 R_{h,N}(t)=\mathcal{O}(th^{N-1}),\ \ \ \  S_{h,N}(t)=\mathcal{O}(th^{N-1}),
\label{remainder}%
\end{equation}
and the coefficient functions as well as all their derivatives are
bounded by
\begin{equation}%
\begin{array}
[c]{rll}%
&\zeta_0(t)=\mathcal{O} (1  ),\quad  &\eta_0(t)= \mathcal{O}(1  ), \\
&\zeta_j(t)=\mathcal{O}\Big(\frac{h^2\bar{b}_1(h\omega_j)}{\sin^2(\frac{1}{2}h\omega_j)
}\Big)=\mathcal{O}(h),\quad
&\eta_j(t)=\mathcal{O}\Big(\frac{h^2b_1(h\omega_j)}{\sin^2(\frac{1}{2}h\omega_j)
}\Big)=\mathcal{O}(h),  \\
&\dot{\zeta}_j^{\pm\langle
j\rangle}(t)=\mathcal{O}\Big(\frac{h^2\bar{b}_1(h\omega_j)}{\sin(\pm
\omega_j h)}\Big)=\mathcal{O}(h^{3/2}), \quad  &
\dot{\eta}_j^{\pm\langle
j\rangle}(t)=\mathcal{O}\Big(\frac{h^2b_1(h\omega_j)}{\sin(\pm
\omega_j
h)}\Big)=\mathcal{O}(h^{3/2}),\\
&\zeta_j^{\pm\langle j\rangle}(t)=\mathcal{O}(\epsilon), \quad  &
\eta_j^{\pm\langle
j\rangle}(t)=\mathcal{O}(1),\\
&\zeta_0^k(t)=\mathcal{O}\Big(h\epsilon^{|k|}\Big),\quad
&\eta_0^k(t)=\mathcal{O}\Big(h\epsilon^{|k|}\Big),\ \ k\in\mathcal{N}^*, \\
&\zeta_j^k(t)=\mathcal{O}\Big( h\epsilon^{|k|}\bar{b}_1(h\omega_j)
\Big)=\mathcal{O}\Big(h\epsilon^{|k|}\Big),\quad
&\eta_j^k(t)=\mathcal{O}\Big(h\epsilon^{|k|}b_1(h\omega_j)
\Big)=\mathcal{O}\Big(h\epsilon^{|k|}\Big),\\
&\ & \qquad \qquad \qquad \qquad \qquad \qquad k\neq\pm\langle
j\rangle,
\end{array}
\label{coefficient func}%
\end{equation}
for $j=1,\ldots,l$. Moreover, we have
$\zeta^{-k}=\overline{\zeta^{k}}$ and
$\eta^{-k}=\overline{\eta^{k}}$. The constants symbolised by the
notation are independent of $h$ and $\omega$, but depend on $E,\ N,\
c_0$ and $T$.
\end{theo}
\emph{Proof.} We will prove that there exist  two functions
\begin{equation}
\begin{aligned} &q_{h}(t)= \zeta(t)+\sum\limits_{k\in\mathcal{N}^*} \mathrm{e}^{\mathrm{i}(k \cdot \omega) t}\zeta^k(t),\
 \ \ \ p_{h}(t)= \eta(t)+\sum\limits_{k\in\mathcal{N}^*} \mathrm{e}^{\mathrm{i}(k \cdot \omega) t}\eta^k(t)
\end{aligned}
\label{MFE-1}%
\end{equation}
with smooth  coefficients $\zeta,\ \zeta^k,\ \eta,\ \eta^k$,  such
that, for $t=nh$,
\begin{equation*}
\begin{aligned} &q^{n}=q_{h}(t)+\mathcal{O}(h^{N}),\ \ \ \  p^{n}= p_{h}(t)+\mathcal{O}(h^{N}).\\
\end{aligned}
\end{equation*}

\vskip2mm \textbf{Construction of the coefficients functions.}

 $\bullet$  For the first   term  of     \eqref{methods}, we look for the  function
\begin{equation}
\begin{aligned} &q^{n+\frac{1}{2}}:=\tilde{q}_{h}(t+\frac{h}{2})=\xi(t+\frac{h}{2})+\sum\limits_{k\in\mathcal{N}^*}
\mathrm{e}^{\mathrm{i}(k \cdot \omega) t}\xi^k(t+\frac{h}{2})
\end{aligned}
\label{MFE-2}%
\end{equation}
for $Q^{n+\frac{1}{2}}$ in the numerical integrator \eqref{methods}.
Inserting  \eqref{MFE-1} and \eqref{MFE-2} into the first  term
 of \eqref{methods}  and comparing the coefficients of $\mathrm{e}^{\mathrm{i}(k \cdot
\omega) t}$, one gets
\begin{equation*}
\begin{aligned} &\xi(t+\frac{h}{2})=\cos(\frac{1}{2}h\Omega)\zeta(t)+\frac{1}{2}h
\mathrm{sinc}(\frac{1}{2}h\Omega)\eta(t),\\
&\xi^k(t+\frac{h}{2})=
\cos(\frac{1}{2}h\Omega)\zeta^k(t)+\frac{1}{2}h
\mathrm{sinc}(\frac{1}{2}h\Omega)\eta^k(t).
\end{aligned}
\end{equation*}

 $\bullet$ For the second term of     \eqref{methods},
by the symmetry of the integrator, we obtain
\begin{equation}
\begin{aligned}&q^{n+1}-2\cos(h\Omega)q^{n}+q^{n-1}=h^2\bar{b}_1(h\Omega)\big[g(q^{n+\frac{1}{2}})+g(q^{n-\frac{1}{2}})\big],\\
\end{aligned}\label{MFE-q2}%
\end{equation}
where  $q^{n-\frac{1}{2}}$ is defined by
$q^{n-\frac{1}{2}}:=\tilde{q}_{h}(t-\frac{h}{2})=\xi(t-\frac{h}{2})+\sum\limits_{k\in\mathcal{N}^*}
\mathrm{e}^{\mathrm{i}(k \cdot \omega)
 t}\xi^k(t-\frac{h}{2})
$ with
\begin{equation*}
\begin{aligned} &\xi(t-\frac{h}{2})=\cos(\frac{1}{2}h\Omega)\zeta(t)-\frac{1}{2}h
\mathrm{sinc}(\frac{1}{2}h\Omega)\eta(t),\\
&\xi^k(t-\frac{h}{2})=
\cos(\frac{1}{2}h\Omega)\zeta^k(t)-\frac{1}{2}h
\mathrm{sinc}(\frac{1}{2}h\Omega)\eta^k(t).
\end{aligned}
\end{equation*}
Inserting the expansions into \eqref{MFE-q2}, we obtain
\begin{equation*}
\begin{aligned}&q_{h}(t+h)-2\cos(h\Omega)q_{h}(t)+q_{h}(t-h)=h^2\bar{b}_1(h\Omega)\big[g(\tilde{q}_{h}(t+\frac{h}{2}))+g(\tilde{q}_{h}(t-\frac{h}{2})\big].
\end{aligned}
\end{equation*}
By the operator $\mathcal{L}(hD)$ and the Taylor series, we can
rewrite the above formula as
\begin{equation*}
\begin{aligned}&\mathcal{L}(hD)q_{h}(t)=h^2\bar{b}_1(h\Omega)\big[g(\tilde{q}_{h}(t+\frac{h}{2}))+g(\tilde{q}_{h}(t-\frac{h}{2})\big]\\
=&h^2\bar{b}_1(h\Omega)\Big[g(\xi(t+\frac{h}{2}))+\sum\limits_{k\in\mathcal{N}^*}\mathrm{e}^{\mathrm{i}(k
\cdot \omega) t}\sum\limits_{s(\alpha)\sim
k}\frac{1}{m!}g^{(m)}(\xi(t+\frac{h}{2}))(\xi(t+\frac{h}{2}))^{\alpha} \\
&+
g(\xi(t-\frac{h}{2}))+\sum\limits_{k\in\mathcal{N}^*}\mathrm{e}^{\mathrm{i}(k
\cdot \omega) t}\sum\limits_{s(\alpha)\sim
k}\frac{1}{m!}g^{(m)}(\xi(t-\frac{h}{2}))(\xi(t-\frac{h}{2}))^{\alpha}\Big],
\end{aligned} %
\end{equation*}
where the sums   are over all $m\geq1$ and over multi-indices
$\alpha=(\alpha_1,\ldots,\alpha_m)$ with $\alpha_j\in
\mathcal{N}^*$,  and the relation $s(\alpha)\sim k$  means
$s(\alpha)- k\in\mathcal{M}.$ Here,
 an abbreviation for the $m$-tuple
$(\xi^{\alpha_1}(t),\ldots,\xi^{\alpha_m}(t))$ is
   denoted by $(\xi(t))^{\alpha}$.

Inserting  the ansatz \eqref{MFE-1} and comparing the coefficients
of $\mathrm{e}^{\mathrm{i}(k \cdot \omega) t}$  yields
\begin{equation*}
\begin{aligned}&\mathcal{L}(hD)\zeta(t)=h^2\bar{b}_1(h\Omega)\Big[g(\xi(t+\frac{h}{2}))+
\sum\limits_{s(\alpha)\sim0}\frac{1}{m!}g^{(m)}(\xi(t+\frac{h}{2}))(\xi(t+\frac{h}{2}))^{\alpha}\\
&\qquad \qquad\qquad+g(\xi(t-\frac{h}{2}))+
\sum\limits_{s(\alpha)\sim0}\frac{1}{m!}g^{(m)}(\xi(t-\frac{h}{2}))(\xi(t-\frac{h}{2}))^{\alpha}\Big],\\
&\mathcal{L}(hD+\mathrm{i}(k \cdot \omega)
h)\zeta^k(t)=h^2\bar{b}_1(h\Omega)\Big[
\sum\limits_{s(\alpha)\sim k}\frac{1}{m!}g^{(m)}(\xi(t+\frac{h}{2}))(\xi(t+\frac{h}{2}))^{\alpha} \\
&\qquad \qquad\qquad\qquad \qquad+
\sum\limits_{s(\alpha)\sim k}\frac{1}{m!}g^{(m)}(\xi(t-\frac{h}{2}))(\xi(t-\frac{h}{2}))^{\alpha}\Big].\\
\end{aligned} %
\end{equation*}
According to the results of $\mathcal{L}(hD)$ and
$\mathcal{L}(hD+\mathrm{i}(k \cdot \omega) h)$ given in Proposition
\ref{lhd pro}, the dominating terms of $\mathcal{L}(hD) \zeta_0(t)$
and $\mathcal{L}(hD) \zeta_j(t)$ are $ h^2D^2\zeta_0(t)$  and $
4\sin^2(h \omega_j/2) \zeta_j(t),$ respectively. Thus, we obtain
\begin{equation}\label{zeta-re}
\begin{aligned}\ddot{\zeta}_0(t)=&\frac{h^2\bar{b}_1(0)}{h^2 }\Big[g(\xi(t+\frac{h}{2}))+
\sum\limits_{s(\alpha)\sim0}\frac{1}{m!}g^{(m)}(\xi(t+\frac{h}{2}))(\xi(t+\frac{h}{2}))^{\alpha}\\
&+g(\xi(t-\frac{h}{2}))+
\sum\limits_{s(\alpha)\sim0}\frac{1}{m!}g^{(m)}(\xi(t-\frac{h}{2}))(\xi(t\frac{h}{2}))^{\alpha}\Big]_0,\\
\zeta_j(t)=&\frac{h^2\bar{b}_1(h\omega_j)}{4\sin^2(\frac{1}{2}h\omega_j)
}\Big[g(\xi(t+\frac{h}{2}))+
\sum\limits_{s(\alpha)\sim0}\frac{1}{m!}g^{(m)}(\xi(t+\frac{h}{2}))(\xi(t+\frac{h}{2}))^{\alpha}\\
&+g(\xi(t-\frac{h}{2}))+
\sum\limits_{s(\alpha)\sim0}\frac{1}{m!}g^{(m)}(\xi(t-\frac{h}{2}))(\xi(t-\frac{h}{2}))^{\alpha}\Big]_0,\
 j=1,\ldots,l.
\end{aligned} %
\end{equation}

Similarly,   the dominating terms of $\mathcal{L}(hD+\mathrm{i}(k
\cdot \omega) h)\zeta_0^k(t)$ for all $k\in \mathcal{N}^*$ are $
\big(2-2\cos((k \cdot \omega) h)\big)\zeta_0^k(t), $ and the
dominating terms of $\mathcal{L}(hD+\mathrm{i}((k \cdot \omega))
h)\zeta_j^{\pm\langle j\rangle}(t)$ for $j=1,\ldots,l$ are $
2\sin(\pm(\langle j\rangle \cdot \omega) h)I(\mathrm{i} h
D)\zeta_j^{\pm\langle j\rangle}(t). $ We also get the dominating
terms of $\mathcal{L}(hD+\mathrm{i}(k \cdot \omega) h)\zeta_j^k(t)$
for $k\neq\pm\langle j\rangle$: $ \big(2\cos((k \cdot \omega)
h)-2\cos(h \omega_j)\big)\zeta_j^k(t). $ Thus, we have
\begin{equation*}
\begin{aligned} \zeta_0^k(t)=&\frac{h^2\bar{b}_1(0)}{2-2\cos((k \cdot
\omega) h)}\Big( \sum\limits_{s(\alpha)\sim
k}\frac{1}{m!}g^{(m)}(\xi(t+\frac{h}{2}))(\xi(t+\frac{h}{2}))^{\alpha}\\
&+\sum\limits_{s(\alpha)\sim
k}\frac{1}{m!}g^{(m)}(\xi(t-\frac{h}{2}))(\xi(t-\frac{h}{2}))^{\alpha}\Big)_0  ,\\
\dot{\zeta}_j^{\pm\langle
j\rangle}(t)=&\frac{h^2\bar{b}_1(h\omega_j)}{2\sin(\pm \omega_j h)}
\Big(\sum\limits_{s(\alpha)\sim
k}\frac{1}{m!}g^{(m)}(\xi(t+\frac{h}{2}))(\xi(t+\frac{h}{2}))^{\alpha}\\
&+\sum\limits_{s(\alpha)\sim
k}\frac{1}{m!}g^{(m)}(\xi(t-\frac{h}{2}))(\xi(t-\frac{h}{2}))^{\alpha}\Big)_j  ,\\
 \zeta_j^k(t)=&\frac{h^2\bar{b}_1(h\omega_j)}{(2\cos((k \cdot \omega)
h)-2\cos(h \omega_j)}\Big( \sum\limits_{s(\alpha)\sim
k}\frac{1}{m!}g^{(m)}(\xi(t+\frac{h}{2}))(\xi(t+\frac{h}{2}))^{\alpha}\\
&+\sum\limits_{s(\alpha)\sim
k}\frac{1}{m!}g^{(m)}(\xi(t-\frac{h}{2}))(\xi(t-\frac{h}{2}))^{\alpha}\Big)_j.
\end{aligned} %
\end{equation*}

 $\bullet$ For the third term of  \eqref{methods}, one arrives
 at
\begin{equation}
\begin{aligned}&p_h(t+h)-2\cos(h\Omega)p_h(t)+p_h(t-h)=hb_1(h\Omega)\big[g(\tilde{q}_{h}(t+\frac{h}{2}))-g(\tilde{q}_{h}(t-\frac{h}{2})\big].
\end{aligned}\label{MFE-q2-newtt}%
\end{equation}
With regard to the coefficient functions $\eta^k(t)$, it is true
that\begin{equation*}
\begin{aligned}
\ddot{\eta}_0(t)=&\frac{hb_1(0)}{h^2 }\Big[g(\xi(t+\frac{h}{2}))+
\sum\limits_{s(\alpha)\sim0}\frac{1}{m!}g^{(m)}(\xi(t+\frac{h}{2}))(\xi(t+\frac{h}{2}))^{\alpha}\\
&-g(\xi(t-\frac{h}{2}))-
\sum\limits_{s(\alpha)\sim0}\frac{1}{m!}g^{(m)}(\xi(t-\frac{h}{2}))(\xi(t-\frac{h}{2}))^{\alpha}\Big]_0,\\
\eta_j(t)=&\frac{hb_1(h\omega_j)}{4\sin^2(\frac{1}{2}h\omega_j)
}\Big[g(\xi(t+\frac{h}{2}))+
\sum\limits_{s(\alpha)\sim0}\frac{1}{m!}g^{(m)}(\xi(t+\frac{h}{2}))(\xi(t+\frac{h}{2}))^{\alpha}\\
&-g(\xi(t-\frac{h}{2}))-
\sum\limits_{s(\alpha)\sim0}\frac{1}{m!}g^{(m)}(\xi(t-\frac{h}{2}))(\xi(t-\frac{h}{2}))^{\alpha}\Big]_j,\\
 \eta_0^k(t)=&\frac{hb_1(0)}{2-2\cos((k \cdot \omega) h)}\Big(
\sum\limits_{s(\alpha)\sim
k}\frac{1}{m!}g^{(m)}(\xi(t+\frac{h}{2}))(\xi(t+\frac{h}{2}))^{\alpha}\\
&-\sum\limits_{s(\alpha)\sim
k}\frac{1}{m!}g^{(m)}(\xi(t-\frac{h}{2}))(\xi(t-\frac{h}{2}))^{\alpha}\Big)_0,
\end{aligned} %
\end{equation*}
and
\begin{equation*}
\begin{aligned}
\dot{\eta}_j^{\pm\langle
j\rangle}(t)=&\frac{hb_1(h\omega_j)}{2\sin(\pm \omega_j h)}
\Big(\sum\limits_{s(\alpha)\sim
k}\frac{1}{m!}g^{(m)}(\xi(t+\frac{h}{2}))(\xi(t+\frac{h}{2}))^{\alpha}\\
&-\sum\limits_{s(\alpha)\sim
k}\frac{1}{m!}g^{(m)}(\xi(t-\frac{h}{2}))(\xi(t-\frac{h}{2}))^{\alpha}\Big)_j ,\\
 \eta_j^k(t)=&\frac{hb_1(h\omega_j)}{(2\cos((k \cdot \omega)
h)-2\cos(h \omega_j)}\Big(\sum\limits_{s(\alpha)\sim
k}\frac{1}{m!}g^{(m)}(\xi(t+\frac{h}{2}))(\xi(t+\frac{h}{2}))^{\alpha}\\
&-\sum\limits_{s(\alpha)\sim
k}\frac{1}{m!}g^{(m)}(\xi(t-\frac{h}{2}))(\xi(t-\frac{h}{2}))^{\alpha}\Big)_j.\\
\end{aligned} %
\end{equation*}
We now  obtain the ansatz of all the modulated Fourier functions.
 Since
the series in the ansatz usually diverge, in this paper we
 truncate them after the $\mathcal{O}(h^{N+1})$ terms (see \cite{Hairer00,Hairer16}).

\vskip2mm

\textbf{Initial values.}  It follows from the conditions
$p_h(0)=p^0$ and $q_h(0)=q^0$ that
\begin{equation*}
\begin{aligned}
&p^0 =\eta(0)+\sum\limits_{k\in\mathcal{N}^*}\eta^k(0) +\mathcal{O}(h^{N}),\ \ q^0=\zeta(0)+\sum\limits_{k\in\mathcal{N}^*}  \zeta^k(0)+\mathcal{O}(h^{N}).\end{aligned} %
\end{equation*}
This implies
\begin{equation}\label{Initial values-1}%
\begin{aligned}
&p_0^0 =\eta_0(0)+\mathcal{O}(h  ),\ \ \ \ q_0^0=\zeta_0(0)+\mathcal{O}(h  ),\\\end{aligned} %
\end{equation}
and
\begin{equation*}
\begin{aligned}
&p_j^0 =\eta_j(0)+\eta^{\pm\langle j\rangle}_j(0) +\mathcal{O}(h
\epsilon),\  \ q_j^0=\zeta_j(0)+\zeta^{\pm\langle
j\rangle}_j(0)+\mathcal{O}(h \epsilon),\end{aligned} %
\end{equation*}
which lead to
\begin{equation}\label{Initial values-2}%
\begin{aligned}
&2\mathrm{Re}(\eta^{ \langle j\rangle}_j(0) )=p_j^0
-\eta_j(0)+\mathcal{O}(h \epsilon),\ \  2\mathrm{Re}(\zeta^{ \langle
j\rangle}_j(0))=q_j^0-\zeta_j(0)+\mathcal{O}(h \epsilon).\end{aligned} %
\end{equation}
Furthermore,  we note that it holds that $p_h(h)=p^1$ and
$q_h(h)=q^1.$ Using the integrator \eqref{methods}, we have
\begin{equation*}\begin{aligned} &q_0^{1}- q_0^{0}=h
 p_0^{0}+h^2\bar{b}_1(0)\big(g(q^{\frac{1}{2}})\big)_0,\ \ p_0^{1}- p_0^{0}= h b_1(0)\big(g(q^{\frac{1}{2}})\big)_0.
 \end{aligned}\end{equation*}
Hence, we arrive at
\begin{equation}\begin{aligned} \label{Initial values-3}
 & \dot{\zeta}_0(0)=
 \eta_0(0)+h\bar{b}_1(0)\big(g(q^{\frac{1}{2}})\big)_0+\mathcal{O}(1),\
 \ \dot{\eta}_0(0)=
b_1(0)\big(g(q^{\frac{1}{2}})\big)_0+\mathcal{O}(1).
 \end{aligned}\end{equation}
The formulae  \eqref{Initial values-1} and \eqref{Initial values-3}
yield the initial values $ \zeta_0(0),\dot{\zeta}_0(0),
\eta_0(0),\dot{\eta}_0(0).$ Therefore, $$\zeta_0(t)=\mathcal{O} (1
),\quad   \eta_0(t)= \mathcal{O}(1  ).$$

Considering again  the integrator \eqref{methods} implies
\begin{equation*} q^{1}-\cos(h\Omega)q^{0}=h
\mathrm{sinc}(h\Omega)p^{0}+h^2\bar{b}_1(h\Omega)g(q^{\frac{1}{2}}).\end{equation*}
On the other hand,  a calculation gives
\begin{equation}\label{Initial values-f1}
\begin{aligned}
&q_j^{1}-\cos(h\omega_j)q_j^{0}=q_j(h)-\cos(h\omega_j)q_j(0)
=\zeta_j(h)+\sum\limits_{k\in\mathcal{N}^*} \mathrm{e}^{\mathrm{i}k
\cdot \omega
h}\zeta_j^k(h)\\
&-\cos(h\omega_j)\Big(\zeta_j(0)+\sum\limits_{k\in\mathcal{N}^*}
\zeta_j^k(0)\Big) =\zeta_j(h)+ \mathrm{e}^{\mathrm{i}  \omega_j
h}\zeta^{ \langle j\rangle}_j(h)+ \mathrm{e}^{-\mathrm{i}  \omega_j
h}\zeta^{ -\langle j\rangle}_j(h)\\
&-\cos(h\omega_j)\Big(\zeta_j(0)+ \zeta_j^{  \langle
j\rangle}(0)+\zeta_j^{ -\langle j\rangle}(0)\Big)+\mathcal{O}(h
\epsilon),
\end{aligned}
\end{equation}
which leads to
\begin{equation*}
\begin{aligned}
&q_j^{1}-\cos(h\omega_j)q_j^{0}\\
=&(1-\cos(h\omega_j))\zeta_j(0)+\mathrm{i}\sin(h\omega_j) (\zeta_j^{
\langle j\rangle}(0)-\zeta_j ^{ -\langle
j\rangle}(0))+\mathcal{O}(h^2)+\mathcal{O}(h \epsilon)
\end{aligned}
\end{equation*}
by expanding the functions $\zeta_j(h),\ \zeta_j^{ \langle
j\rangle}(h)$ and $\zeta_j^{ -\langle j\rangle}(h)$ at $h=0$. From
the fact that
$1-\cos(h\omega_j)=\frac{1}{2}h^2\omega_j^2\mathrm{sinc}^2(h\omega_j/2),$
it follows that
\begin{equation*}
\begin{aligned}
&(1-\cos(h\omega_j))\zeta_j(0)=\frac{1}{2}h^2\omega_j^2\mathrm{sinc}^2(h\omega_j/2))\zeta_j(0)=2\sin^2(h\omega_j/2)\zeta_j(0).
\end{aligned}
\end{equation*}
In the light of  the second formula of \eqref{zeta-re}, we get
another expression of the above result
\begin{equation*}
\begin{aligned}
&(1-\cos(h\omega_j))\zeta_j(0) = \frac{1}{2}
h^2\bar{b}_1(h\omega_j)\Big[g(\xi(\frac{h}{2}))+
\sum\limits_{s(\alpha)\sim0}\frac{1}{m!}g^{(m)}(\xi(\frac{h}{2}))(\xi(\frac{h}{2}))^{\alpha}\\
&+g(\xi(-\frac{h}{2}))+
\sum\limits_{s(\alpha)\sim0}\frac{1}{m!}g^{(m)}(\xi(-\frac{h}{2}))(\xi(-\frac{h}{2}))^{\alpha}\Big]_j.
\end{aligned}
\end{equation*}
Then \eqref{Initial values-f1} has the following form
\begin{equation*}
\begin{aligned}
&\mathrm{i}\sin(h\omega_j) (\zeta_j^{ \langle j\rangle}(0)-\zeta_j
^{ -\langle j\rangle}(0)) =h
\mathrm{sinc}(h\omega_j)p_j^{0}+h^2\bar{b}_1(h\omega_j)(g(q^{\frac{1}{2}}))_j\\&+\frac{1}{2}h^2\bar{b}_1(h\omega_j)\Big[g(\xi(\frac{h}{2}))+
\sum\limits_{s(\alpha)\sim0}\frac{1}{m!}g^{(m)}(\xi(\frac{h}{2}))(\xi(\frac{h}{2}))^{\alpha}\\
&+g(\xi(-\frac{h}{2}))+
\sum\limits_{s(\alpha)\sim0}\frac{1}{m!}g^{(m)}(\xi(-\frac{h}{2}))(\xi(-\frac{h}{2}))^{\alpha}\Big]_j
+\mathcal{O}(h \epsilon),
\end{aligned}
\end{equation*}
which yields
\begin{equation}\label{Initial values-4}2\mathrm{Im}(\zeta_j^{ \langle j\rangle}(0))=\omega_j^{-1}p_j^{0}+\mathcal{O}(\epsilon).\end{equation}
Similarly, it can be obtained that
\begin{equation}\label{Initial values-5}2\mathrm{Im}(\eta_j^{ \langle j\rangle}(0))=-\omega_j q_j^{0}
 +\mathcal{O}(\epsilon).\end{equation} The conditions \eqref{Initial values-2},
\eqref{Initial values-4} and \eqref{Initial values-5} present the
desired initial values $\zeta_j^{ \pm\langle j\rangle}(0)$ and $
\eta_j^{\pm \langle j\rangle}(0)$. This analysis implies
\begin{equation*}
\begin{aligned}
& \zeta_j^{ \pm\langle j\rangle}(t)=\mathcal{O}(  \epsilon),\ \  \eta_j^{ \pm\langle j\rangle}(t)=\mathcal{O}(  1).\\
\end{aligned}
\end{equation*}

\vskip2mm

 \textbf{Bounds.}
Based on the ansatz, the initial values and Assumption \ref{ass}, it
is easy to get the bounds \eqref{coefficient func} of modulated
Fourier functions.

 \textbf{Defect.}  The defect \eqref{remainder} can be obtained  by using
the Lipschitz continuous of the nonlinearity $g$, a discrete
Gronwall lemma and  the standard  convergence estimates (see
\cite{Hairer00,17-new}  and Chap. XIII of \cite{hairer2006} for more
details).
\\

We then complete the proof of this theorem.
 \hfill  $\square$
\section{Almost-invariants of the integrators} \label{sec:Almost-invariants}
In this section, we show that the ERKN integrators have two
almost-invariants.
\subsection{The first almost-invariant}
 Let $\vec{\zeta}=\big(\zeta^{k}\big)_{k\in
\mathcal{N}}$ and $ \vec{\eta}=\big(\eta^{k}\big)_{k\in
\mathcal{N}}.$ The first almost-invariant is given as follows.
\begin{theo}\label{first invariant thm}
Under the conditions of Theorem \ref{energy thm}, there exists a
function $\widehat{\mathcal{H}}[\vec{\zeta},\vec{\eta}]$ such that
\begin{equation*}
\widehat{\mathcal{H}}[\vec{\zeta},\vec{\eta}](t)=\widehat{\mathcal{H}}[\vec{\zeta},\vec{\eta}](0)+\mathcal{O}(th^{N})
\end{equation*}
for $0\leq t\leq T.$ Moreover, $\widehat{\mathcal{H}}$ can be
expressed in
\begin{equation*}\begin{aligned}
\widehat{\mathcal{H}}[\vec{\zeta},\vec{\eta}]=&\frac{1}{4\bar{b}_1}
\eta_0
 ^\intercal\eta_0
+\sum\limits_{j=1}^l2\omega_j^2\mathrm{sinc}(
h\omega_j)\frac{\cos(\frac{1}{2}h\omega_j)}{2\bar{b}_1(h\omega_j)}\big(\zeta_j^{-\langle
j\rangle}\big)^\intercal\zeta_j^{\langle j\rangle}\\
&+\sum\limits_{j=1}^l2h^2\omega_j^2\mathrm{sinc}(
h\omega_j)\frac{\frac{1}{2}\mathrm{sinc}(\frac{1}{2}h\omega_j)}{2b_1(h\omega_j)}\big(\eta_j^{-\langle
j\rangle}\big)^\intercal\eta_j^{\langle j\rangle}
+U(\xi(t))+\mathcal{O}(h).
\end{aligned}
\end{equation*}
\end{theo}
\emph{Proof.} It follows from the proof of Theorem \ref{energy thm}
that
\begin{equation*}
\begin{aligned} &\tilde{q}_{h}(t+\frac{h}{2})=\cos(\frac{1}{2}h\Omega)q_{h}(t)+\frac{1}{2}h \mathrm{sinc}(\frac{1}{2}h\Omega)p_h(t),\\
&\tilde{q}_{h}(t-\frac{h}{2})=\cos(\frac{1}{2}h\Omega)q_{h}(t)-\frac{1}{2}h \mathrm{sinc}(\frac{1}{2}h\Omega)p_h(t),\\
& \mathcal{L}(hD) q_{h}(t)=h^2\bar{b}_1(h\Omega)\big(g(\tilde{q}_{h}(t+\frac{h}{2}))+g(\tilde{q}_{h}(t-\frac{h}{2}))\big)+\mathcal{O}(h^{N}),\\
& \mathcal{L}(hD) p_{h}(t)=hb_1(h\Omega)\big(g(\tilde{q}_{h}(t+\frac{h}{2}))-g(\tilde{q}_{h}(t-\frac{h}{2}))\big)+\mathcal{O}(h^{N}),\\
\end{aligned}
\end{equation*}
where   the following denotations are used:
\begin{equation*}
\begin{aligned}q_{h}(t)=\sum\limits_{k\in\mathcal{N}}q^k_{h}(t),\  \ p_{h}(t)=\sum\limits_{
k\in\mathcal{N}}p^k_{h}(t),\ \
\tilde{q}_{h}(t\pm\frac{1}{2}h)=\sum\limits_{
k\in\mathcal{N}}\tilde{q}^k_{h}(t\pm\frac{1}{2}h)
\end{aligned}
\end{equation*}
with
\begin{equation*}
\begin{aligned}q^k_{h}(t)=\mathrm{e}^{\mathrm{i}(k \cdot \omega) t}\zeta^k(t),\ \ p^k_{h}(t)= \mathrm{e}^{\mathrm{i}(k \cdot \omega)
t}\eta^k(t),\ \
\tilde{q}^k_{h}(t\pm\frac{1}{2}h)=\mathrm{e}^{\mathrm{i}(k \cdot
\omega) t}\xi^k(t\pm\frac{1}{2}h).
\end{aligned}
\end{equation*}
This yields
\begin{equation}
\begin{aligned} &\tilde{q}^k_{h}(t+\frac{h}{2})=\cos(\frac{1}{2}h\Omega)q^k_{h}(t)+\frac{1}{2}h \mathrm{sinc}(\frac{1}{2}h\Omega)p^k_h(t),\\
&\tilde{q}^k_{h}(t-\frac{h}{2})=\cos(\frac{1}{2}h\Omega)q^k_{h}(t)-\frac{1}{2}h \mathrm{sinc}(\frac{1}{2}h\Omega)p^k_h(t),\\
& \mathcal{L}(hD)
q^k_{h}(t)=-h^2\bar{b}_1(h\Omega)\Big(\nabla_{q^{-k}}\mathcal{U}(\tilde{q}
(t+\frac{h}{2}))
+\nabla_{q^{-k}}\mathcal{U}(\tilde{q} (t-\frac{h}{2}))\Big)+\mathcal{O}(h^{N}),\\
& \mathcal{L}(hD)
p^k_{h}(t)=-hb_1(h\Omega)\Big(\nabla_{q^{-k}}\mathcal{U}(\tilde{q}
(t+\frac{h}{2})) -\nabla_{q^{-k}}\mathcal{U}(\tilde{q}
(t-\frac{h}{2}))\Big)+\mathcal{O}(h^{N}),
\end{aligned}
\label{methods-inva-nnew}%
\end{equation}
where $\mathcal{U}(\tilde{q})$ is defined as
\begin{equation}
\mathcal{U}(\tilde{q}(t\pm\frac{h}{2}))=U(\tilde{q}_h(t\pm\frac{h}{2}))+
\sum\limits_{s(\alpha)\sim0}\frac{1}{m!}U^{(m)}(\tilde{q}_h(t\pm\frac{h}{2}))
(\tilde{q}_h(t\pm\frac{h}{2}))^{\alpha}
\label{newuu}%
\end{equation}
with $
\tilde{q}(t\pm\frac{h}{2})=\big(\tilde{q}^{k}_h(t\pm\frac{h}{2})\big)_{k\in
\mathcal{N}}.$ Thence, the following result is obtained
\begin{equation*}
\begin{aligned}
&\frac{1}{2}\frac{d}{dt}
\Big(\mathcal{U}(\tilde{q}(t+\frac{h}{2}))+\mathcal{U}(\tilde{q}(t-\frac{h}{2}))\Big)\\
=&\frac{1}{2}\sum\limits_{k\in\mathcal{N}}\Big[
\big(\dot{\tilde{q}}^{-k}_{h}(t+\frac{h}{2})\big)^\intercal
\nabla_{q^{-k}}\mathcal{U}(\tilde{q}(t+\frac{h}{2}))
+\big(\dot{\tilde{q}}^{-k}_{h}(t-\frac{h}{2})\big)^\intercal
\nabla_{q^{-k}}\mathcal{U}(\tilde{q}(t-\frac{h}{2}))\Big]\\
=&\frac{1}{2}\sum\limits_{k\in\mathcal{N}}\Big[
\big(\cos(\frac{1}{2}h\Omega)\dot{q}^{-k}_{h}(t)+\frac{1}{2}h
\mathrm{sinc}(\frac{1}{2}h\Omega)\dot{p}^{-k}_h(t)\big)^\intercal
\nabla_{q^{-k}}\mathcal{U}(\tilde{q}(t+\frac{h}{2}))\\
& +\big(\cos(\frac{1}{2}h\Omega)\dot{q}^{-k}_{h}(t)-\frac{1}{2}h
\mathrm{sinc}(\frac{1}{2}h\Omega)\dot{p}^{-k}_h(t)\big)^\intercal
\nabla_{q^{-k}}\mathcal{U}(\tilde{q}(t-\frac{h}{2}))\Big].
\end{aligned}
\end{equation*}
 By the last two equations of \eqref{methods-inva-nnew}, this formula
 becomes
\begin{equation*}
\begin{aligned}
&\frac{1}{2}\frac{d}{dt}
\Big(\mathcal{U}(\tilde{q}(t+\frac{h}{2}))+\mathcal{U}(\tilde{q}(t-\frac{h}{2}))\Big)
=\frac{1}{2}\sum\limits_{k\in\mathcal{N}}\Big[
\big(\dot{q}^{-k}_{h}(t)\big)^\intercal
\cos(\frac{1}{2}h\Omega)(-h^2\bar{b}_1(h\Omega))^{-1}\\
&\mathcal{L}(hD) q^k_{h}(t)+ \big(\dot{p}^{-k}_h(t)\big)^\intercal
\frac{1}{2}h \mathrm{sinc}(\frac{1}{2}h\Omega)(-h
b_1(h\Omega))^{-1}\mathcal{L}(hD)
p^k_{h}(t)\Big]+\mathcal{O}(h^{N}).
\end{aligned}
\end{equation*}
Rewrite it in the quantities $\zeta_h^k(t),\ \eta_h^k(t)$
\begin{equation}
\begin{aligned}
&\frac{1}{2}\frac{d}{dt}
\Big(\mathcal{U}(\xi_h(t+\frac{h}{2}))+\mathcal{U}(\xi_h(t-\frac{h}{2}))\Big)
+\frac{1}{2}\sum\limits_{k\in\mathcal{N}}\Big[
\big(\dot{\zeta}_h^{-k} (t)- \textmd{i}(k\cdot\omega) \zeta_h^{-k}
(t)\big)^\intercal\\&
\cos(\frac{1}{2}h\Omega)(h^2\bar{b}_1(h\Omega))^{-1}\mathcal{L}(hD+\mathrm{i}(k
\cdot \omega) h) \zeta^k_{h}(t) + \big(\dot{\eta}_h^{-k} (t)-
\textmd{i}(k\cdot\omega) \eta_h^{-k} (t)\big)^\intercal
\\
&\frac{1}{2}h \mathrm{sinc}(\frac{1}{2}h\Omega)(h
b_1(h\Omega))^{-1}\mathcal{L}(hD+\mathrm{i}(k \cdot \omega) h)
\eta^k_{h}(t)\Big] =\mathcal{O}(h^{N}),
\end{aligned}
\label{duu}%
\end{equation}
where $
\xi_h(t\pm\frac{h}{2})=\big(\xi_h^{k}(t\pm\frac{h}{2})\big)_{k\in
\mathcal{N}}.$

With the analysis given in Section XIII of \cite{hairer2006},  it is
known that the left-hand side of \eqref{duu} is a total derivative
and its  construction is given as follows. According to the ``magic
formulas" on p. 508 of \cite{hairer2006} and the bounds of Theorem
\ref{energy thm}, we have
\begin{equation*}
\begin{aligned}\widehat{\mathcal{H}}[\vec{\zeta},\vec{\eta}]=&\frac{1}{4\bar{b}_1} \dot{\zeta}_0
 ^\intercal\dot{\zeta}_0
+\sum\limits_{j=1}^l2\omega_j^2\mathrm{sinc}(
h\omega_j)\frac{\cos(\frac{1}{2}h\omega_j)}{2\bar{b}_1(h\omega_j)}\big(\zeta_j^{-\langle
j\rangle}\big)^\intercal\zeta_j^{\langle j\rangle}\\
&+\sum\limits_{j=1}^l2h^2\omega_j^2\mathrm{sinc}(
h\omega_j)\frac{\frac{1}{2}\mathrm{sinc}(\frac{1}{2}h\omega_j)}{2b_1(h\omega_j)}\big(\eta_j^{-\langle
j\rangle}\big)^\intercal\eta_j^{\langle j\rangle}\\
&
+\frac{1}{2}\big(U(\xi(t+\frac{1}{2}h))+U(\xi(t-\frac{1}{2}h))\big)+\mathcal{O}(h)\\
=&\frac{1}{4\bar{b}_1}  \eta_{h,1}
 ^\intercal \eta_{h,1}
+\sum\limits_{j=1}^l2\omega_j^2\mathrm{sinc}(
h\omega_j)\frac{\cos(\frac{1}{2}h\omega_j)}{2\bar{b}_1(h\omega_j)}\big(\zeta_j^{-\langle
j\rangle}\big)^\intercal\zeta_j^{\langle j\rangle}\\
&+\sum\limits_{j=1}^l2h^2\omega_j^2\mathrm{sinc}(
h\omega_j)\frac{\frac{1}{2}\mathrm{sinc}(\frac{1}{2}h\omega_j)}{2b_1(h\omega_j)}\big(\eta_j^{-\langle
j\rangle}\big)^\intercal\eta_j^{\langle j\rangle}
+U(\xi(t))+\mathcal{O}(h),
\end{aligned}
\end{equation*}
where the fact that $\dot{\zeta}_{h,1}=\eta_{h,1}+\mathcal{O}(h)$
  is used. The proof is complete.
 \hfill  $\square$

\subsection{The second  almost-invariant}

For $\mu\in \mathbb{R}^l$ and
$\tilde{q}(t\pm\frac{1}{2}h)=(\tilde{q}^k_{h}(t\pm\frac{1}{2}h))_{k\in\mathcal{N}}$,
let
 $$S_{\mu}(\tau)\tilde{q}(t\pm\frac{1}{2}h)=(\mathrm{e}^{\mathrm{i}(k
\cdot \mu)
\tau}\tilde{q}^k_{h}(t\pm\frac{1}{2}h))_{k\in\mathcal{N}},\ \ \
\tau\in \mathbb{R}.$$ Inserting
$S_{\mu}(\tau)\tilde{q}(t\pm\frac{1}{2}h)$ into  \eqref{newuu}
yields
\begin{equation}
\begin{aligned}
&\mathcal{U}(S_{\mu}(\tau)\tilde{q}(t\pm\frac{1}{2}h))=
U(\tilde{q}_h (t\pm\frac{1}{2}h))
+\sum\limits_{s(\alpha)\sim0}\frac{\mathrm{e}^{\mathrm{i}(s(\alpha)
\cdot \mu ) \tau} }{m!}U^{(m)}(\tilde{q}_h
(t\pm\frac{1}{2}h))\\&(\mathrm{e}^{\mathrm{i}(\alpha_1 \cdot \mu)
\tau}\tilde{q}_h^{\alpha_1}(t\pm\frac{1}{2}h),\ldots,\mathrm{e}^{\mathrm{i}(\alpha_m
\cdot \mu)  \tau}\tilde{q}_h^{\alpha_m}(t\pm\frac{1}{2}h)).
\end{aligned}
\label{dnewuu}%
\end{equation}

If $\mu\perp \mathcal{M}$, then it follows from the relation
$s(\alpha)\sim0$ that  $s(\alpha)\cdot \mu =0$. This means that  the
expression \eqref{dnewuu} is independent of $\tau$. Therefore, we
have
$$0=\frac{d}{d\tau}\mid_{\tau=0}\mathcal{U}(S_{\mu}(\tau)\tilde{q}(t\pm\frac{1}{2}h))=\sum\limits_{k\in\mathcal{N}}\mathrm{i}(k \cdot
\mu) (\tilde{q}_h^{k}(t\pm\frac{1}{2}h))^\intercal
\nabla_{q^{k}}\mathcal{U}( \tilde{q} (t\pm\frac{1}{2}h)).$$

If $\mu$ is not orthogonal to $ \mathcal{M}$, this means that some
terms in the sum of \eqref{dnewuu}  depend on $\tau$.  For these
terms with $s(\alpha)\in \mathcal{M}$ and $s(\alpha)\cdot \mu
\neq0$, we have $| s(\alpha)|\geq M=\min\{| k|:0\neq k\in
\mathcal{M}\}$ and if $\mu\perp \mathcal{M}_N:=\{k\in\mathcal{M}:\
|k|\leq N\}$, then $|s(\alpha)|\geq N+1$. This result as well as the
bounds \eqref{coefficient func} then implies
\begin{equation*}
\begin{aligned}
 \sum\limits_{k\in\mathcal{N}}\mathrm{i}(k \cdot
\mu)(\tilde{q}_h^{k}(t\pm\frac{1}{2}h))^\intercal
\nabla_{q^{k}}\mathcal{U}( \tilde{q} (t\pm\frac{1}{2}h)) =\left\{
\begin{aligned}
 &\mathcal{O}(\epsilon^M),\ \ \ \ \textmd{for arbitrary}\ \mu,\\
 &\mathcal{O}(\epsilon^{N+1}), \  \textmd{for}\ \mu\perp \mathcal{M}_N.\end{aligned}\right.
\end{aligned}
\end{equation*}
Therefore, we obtain
\begin{equation*}
\begin{aligned}&\mathcal{O}(h^{N})+\mathcal{O}(\epsilon^{M-1})=\frac{\mathrm{i}}{\epsilon}\Big[\frac{d}{d\tau}\mid_{\tau=0}\mathcal{U}(S_{\mu}(\tau)\tilde{q}(t+\frac{1}{2}h))
+\mathcal{U}(S_{\mu}(\tau)\tilde{q}(t-\frac{1}{2}h))\Big]\\
=&\frac{\mathrm{i}}{\epsilon}\sum\limits_{k\in\mathcal{N}}(k \cdot
\mu) \Big[ (\tilde{q}_h^{k}(t+\frac{1}{2}h))^\intercal
\nabla_{q^{k}}\mathcal{U}( \tilde{q}
(t+\frac{1}{2}h))+(\tilde{q}_h^{k}(t-\frac{1}{2}h))^\intercal
\nabla_{q^{k}}\mathcal{U}( \tilde{q} (t-\frac{1}{2}h))
\Big],\end{aligned}
\end{equation*} and the $\mathcal{O}(\epsilon^{M-1})$ term
can be removed for $\mu\perp \mathcal{M}_N.$

In a similar way to the proof of Theorem \ref{first invariant thm},
the above analysis yields the following second almost-invariant.
\begin{theo}\label{second invariant thm}
Under the conditions of Theorem \ref{first invariant thm}, there
exists a function $\widehat{\mathcal{I}}[\vec{\zeta},\vec{\eta}]$
such that
\begin{equation*}
\widehat{\mathcal{I}}_{\mu}[\vec{\zeta},\vec{\eta}](t)=\widehat{\mathcal{I}}_{\mu}[\vec{\zeta},\vec{\eta}](0)+\mathcal{O}(th^{N})+\mathcal{O}(t\epsilon^{M-1})
\end{equation*}
for all $\mu\in \mathbb{R}^l$ and $0\leq t\leq T.$ They satisfy
\begin{equation*}
\widehat{\mathcal{I}}_{\mu}[\vec{\zeta},\vec{\eta}](t)=\widehat{\mathcal{I}}_{\mu}[\vec{\zeta},\vec{\eta}](0)+\mathcal{O}(th^{N})
\end{equation*}
 for $ \mu\perp \mathcal{M}_N$ and $0\leq t\leq T.$
 Moreover, $\widehat{\mathcal{I}}$ can be expressed in
\begin{equation*}\begin{aligned}
\widehat{\mathcal{I}}_{\mu}[\vec{\zeta},\vec{\eta}]
=&\sum\limits_{j=1}^l2\omega_j^2\mathrm{sinc}(
h\omega_j)\frac{\cos(\frac{1}{2}h\omega_j)}{2\bar{b}_1(h\omega_j)}\frac{\mu_j}{\lambda_j}\big(\zeta_j^{-\langle
j\rangle}\big)^\intercal\zeta_j^{\langle j\rangle}\\
&+\sum\limits_{j=1}^l2h^2\omega_j^2\mathrm{sinc}(
h\omega_j)\frac{\frac{1}{2}\mathrm{sinc}(\frac{1}{2}h\omega_j)}{2b_1(h\omega_j)}\frac{\mu_j}{\lambda_j}\big(\eta_j^{-\langle
j\rangle}\big)^\intercal\eta_j^{\langle j\rangle} +\mathcal{O}(h).
\end{aligned}
\end{equation*}
\end{theo}

\section{Long-time near-conservation of total and oscillatory
energy}\label{sec:Long-time near-conservation} Based on the previous
analysis of this paper and following \cite{Cohen05,17-new} and
Section XIII of \cite{hairer2006}, it is easy to obtain the
following result.

\begin{theo}\label{HHthm} Under the
conditions of Theorem \ref{first invariant thm} and the  additional
 condition
\begin{equation}\label{new-cond}
\begin{aligned}
\mathrm{sinc}( h\Omega)\frac{ \cos(\frac{1}{2}h\Omega)
}{2\bar{b}_1(h\Omega)}+h^2\Omega^2\mathrm{sinc}(
h\Omega)\frac{\frac{1}{2} \mathrm{sinc}(\frac{1}{2}h\Omega)
}{2b_1(h\Omega)}= I,
\end{aligned}
\end{equation}
it holds that
\begin{equation}\label{HIHI}
\begin{aligned}
&\widehat{\mathcal{H}}[\vec{\zeta},\vec{\eta}](nh)=H(q_n,p_n)+\mathcal{O}(h),\\
&\widehat{\mathcal{I}}_{\langle
j\rangle}[\vec{\zeta},\vec{\eta}](nh)=I_{j}(q_n,p_n)+\mathcal{O}(h),
\end{aligned}
\end{equation}
where  the constants symbolized by $\mathcal{O}$  depend on $N,\ T$
and the constants in the assumptions.
\end{theo}

\begin{rem}
It is noted that  the symmetry condition \eqref{sym cond} and the
condition \eqref{new-cond}  determine a  symmetric and symplectic
 ERKN integrator (ERKN3 presented in next section). The appearance \eqref{new-cond}   is obtained
 by  requiring  the  almost-invariants   $\widehat{\mathcal{H}}$ and $\widehat{\mathcal{I}}_{\langle
j\rangle}$ to be close to the energies $H$ and $I_j$, respectively.
The mechanism is not in any obvious way related to symplecticity.
 The same coincidence happens in the analysis of trigonometric integrators  in \cite{Hairer00}.
\end{rem}

The near conservation of $H$  and $I$  over long time intervals is
given by the following theorem.
\begin{theo}\label{Long-time near-conservation thm} Under the
conditions of Theorem \ref{HHthm},   we have
\begin{equation*}
\begin{aligned}
H(q^n,p^n)&=H(q^0,p^0)+\mathcal{O}(h),\\
I_{j}(q^n,p^n)&=I_{j}(q^0,p^0)+\mathcal{O}(h)\\
\end{aligned}
\end{equation*}
for $0\leq nh\leq h^{-N+1}$ and $j=1,2\ldots,l$. The constants
symbolized by $\mathcal{O}$ are independent of $n,\ h,\ \Omega$, but
depend on $N, T$ and the constants in the assumptions.
\end{theo}

To be able to treat the ERKN integrators which are symmetric but do
not satisfy \eqref{new-cond}, we  consider  the modified energies
\begin{equation*}
\begin{aligned}
& H^{*}(q,p)= H(q,p)+\sum\limits_{j=1}^l(\sigma(\xi_j)-1)I_j(q,p)
\end{aligned}
\end{equation*}
and
\begin{equation*}
\begin{aligned}
& I_{\mu}^{*}(q,p)=  \sum\limits_{j=1}^l
\sigma(\xi_j)\frac{\mu_j}{\lambda_j} I_j(q,p),
\end{aligned}
\end{equation*}
 where $\sigma$ is defined by
$$
 \sigma(\xi_j):=\mathrm{sinc}(\xi_j)\frac{ \cos(\frac{1}{2}\xi_j)
}{2\bar{b}_1(\xi_j)}+\xi_j^2\mathrm{sinc}( \xi_j)\frac{\frac{1}{2}
\mathrm{sinc}(\frac{1}{2}\xi_j)
}{2b_1(\xi_j)}=\frac{\cos(\frac{1}{2}\xi_j)}{b_1(\xi_j)}.$$

\noindent We then obtain the following result.
\begin{theo}\label{HHthm-new} Under the
conditions of Theorem \ref{first invariant thm}  and that
$\bar{b}_1=\frac{1}{2}$, it holds that
\begin{equation*}
\begin{aligned}
&\widehat{\mathcal{H}}[\vec{\zeta},\vec{\eta}](nh)=H^{*}(q_n,p_n)+\mathcal{O}(h),\\
&\widehat{\mathcal{I}}_{\mu}[\vec{\zeta},\vec{\eta}](nh)=
I_{\mu}^{*}(q_n,p_n)+\mathcal{O}(h),
\end{aligned}
\end{equation*}
Moreover,    we have
  \begin{equation*}
\begin{aligned}
H^{*}(q^n,p^n)&=H^{*}(q^0,p^0)+\mathcal{O}(h),\\
I_{\mu}^{*}(q^n,p^n)&=I_{\mu}^{*}(q^0,p^0)+\mathcal{O}(h)\\
\end{aligned}
\end{equation*}
for $0\leq nh\leq h^{-N+1}$,   $\mu\in \mathbb{R}^l$ and $ \mu\perp
\mathcal{M}_N$. The constants symbolized by $\mathcal{O}$ are
independent of $n,\ h,\ \Omega$, but depend on $N, T$ and the
constants in the assumptions.
\end{theo}

\section{Numerical examples}\label{sec:examples}

As examples, we present four practical one-stage explicit ERKN
integrators whose coefficients are given in Table \ref{praERKN}.
From Theorem \ref{order conditions}, it follows that all these
integrators are of order two. According to Theorems \ref{symmetric
thm} and \ref{symplectic thm}, the symmetry and symplecticness for
these integrators are shown in Table \ref{praERKN}.

\renewcommand\arraystretch{2.3}
\begin{table}$$
\begin{array}{|c|c|c|c|c|c|c|c|}
\hline
\text{Methods} &c_1  &\bar{b}_1(V)   &b_1(V)   & \text{Symmetric}  &  \text{Symplectic}  \\
\hline
\text{ERKN1} & \frac{1}{2} & \phi_{2}(V)& \phi_{0}(V/4)   & \text{Non}  &\text{Non}   \cr
\text{ERKN2} & \frac{1}{2} & \frac{1}{2} \phi_{0}(V/4) \phi_{1}(V) & \phi^3_{0}(V/4)   & \text{Symmetric} & \text{Non}  \cr
\text{ERKN3} & \frac{1}{2} & \frac{1}{2} \phi_{1}(V/4) & \phi_{0}(V/4)   & \text{Symmetric} & \text{Symplectic} \cr
\text{ERKN4} & \frac{1}{2} &\frac{1}{2} \phi_{1}(V)\phi_{1}(V/4) &\phi_{1}(V)\phi_{0}(V/4)   & \text{Symmetric} & \text{Non} \cr
 \hline
\end{array}
$$
\caption{Four one-stage explicit ERKN integrators.} \label{praERKN}
\end{table}

In order to  illustrate the numerical conservation of
 energies for these four integrators, a Hamiltonian \eqref{H} with
$l=3,\ \lambda=(1,\sqrt{2},2)$ is considered (see \cite{Cohen05}).
It follows from the discussion in \cite{Cohen05} that there is the
$1:2$ resonance between $\lambda_1$ and $\lambda_3$: $\mathcal{M}=
\{(-2k_3,0,k_3):\ k_3\in \mathbb{Z}\}.$ For this problem,  the
dimension of $q_1 = (q_{11}, q_{12})$ is assumed to be 2 and all the
other $q_j$ are assumed to be   1. We choose  $ \epsilon^{-1}=
\omega = 70$, the potential
$$U(q) = (0.001q_0 + q_{11} + q_{22} + q_2 + q_3)^
4,$$ and
$$q(0) = (1,0.3\epsilon,0.8\epsilon,-1.1\epsilon,0.7\epsilon),\ \ p(0) = (-0.75,0.6,0.7,-0.9,0.8)$$
as initial values. For $\lambda=(1,\sqrt{2},2),$ we consider
$$\mu=(1,0,2)\ \ \ \textmd{and} \ \ \ \mu=(0,\sqrt{2},0)$$
 for $I_{\mu}$ and  the corresponding results  are
$$I_{\mu}=I_1+I_3\ \ \ \textmd{and} \ \ \ I_{\mu}=I_2.$$

\begin{figure}[ptb]
\centering\tabcolsep=2mm
\begin{tabular}
[c]{ccc}%
\includegraphics[width=13cm,height=3.5cm]{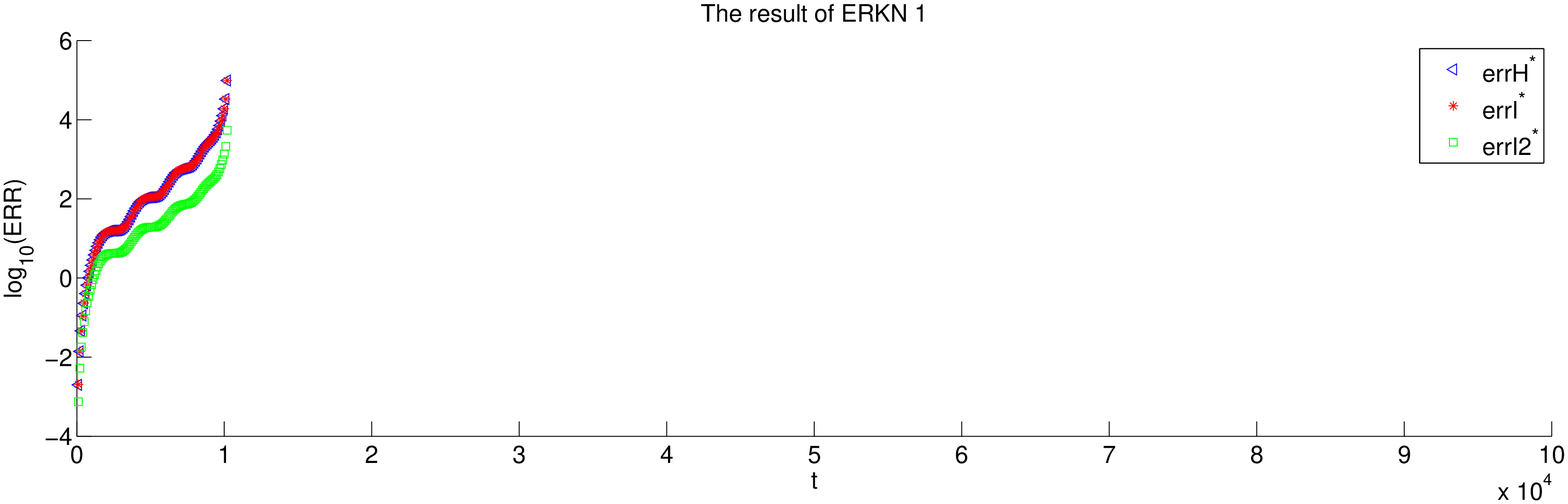}\\
\includegraphics[width=13cm,height=3.5cm]{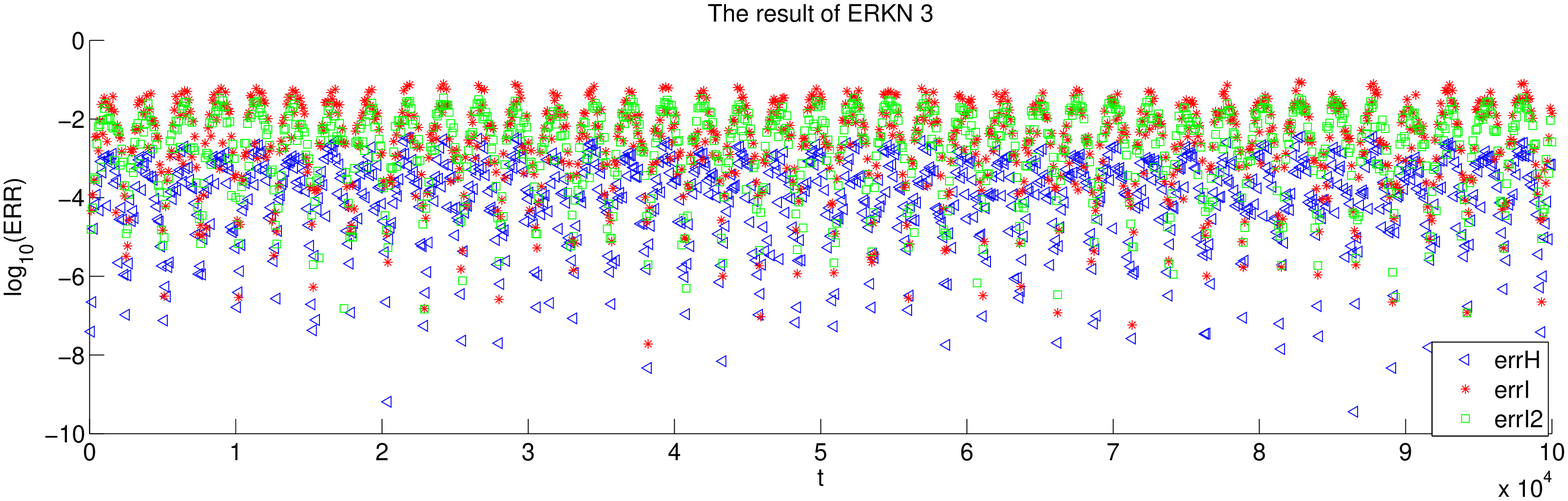}\\
\end{tabular}
\caption{The    errors   against $t$   for ERKN 1 (up) and 3 (down).}%
\label{fig1}%
\end{figure}

\begin{figure}[ptb]
\centering\tabcolsep=2mm
\begin{tabular}
[c]{ccc}%
\includegraphics[width=13cm,height=3.5cm]{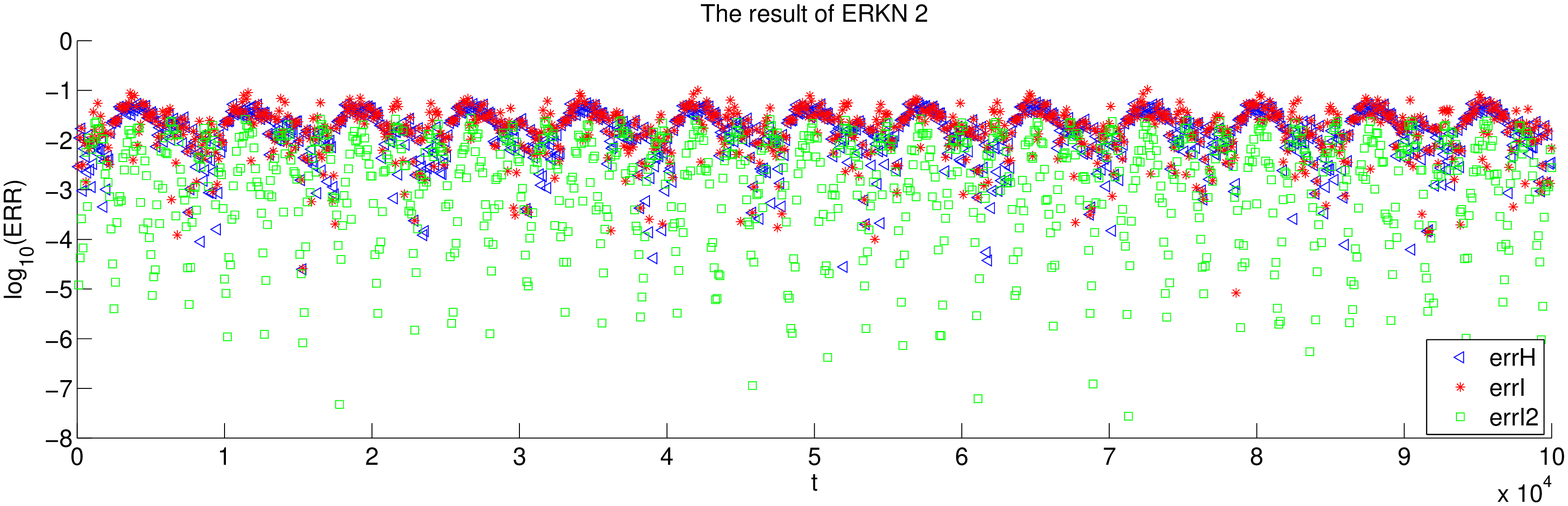}  \\
\includegraphics[width=13cm,height=3.5cm]{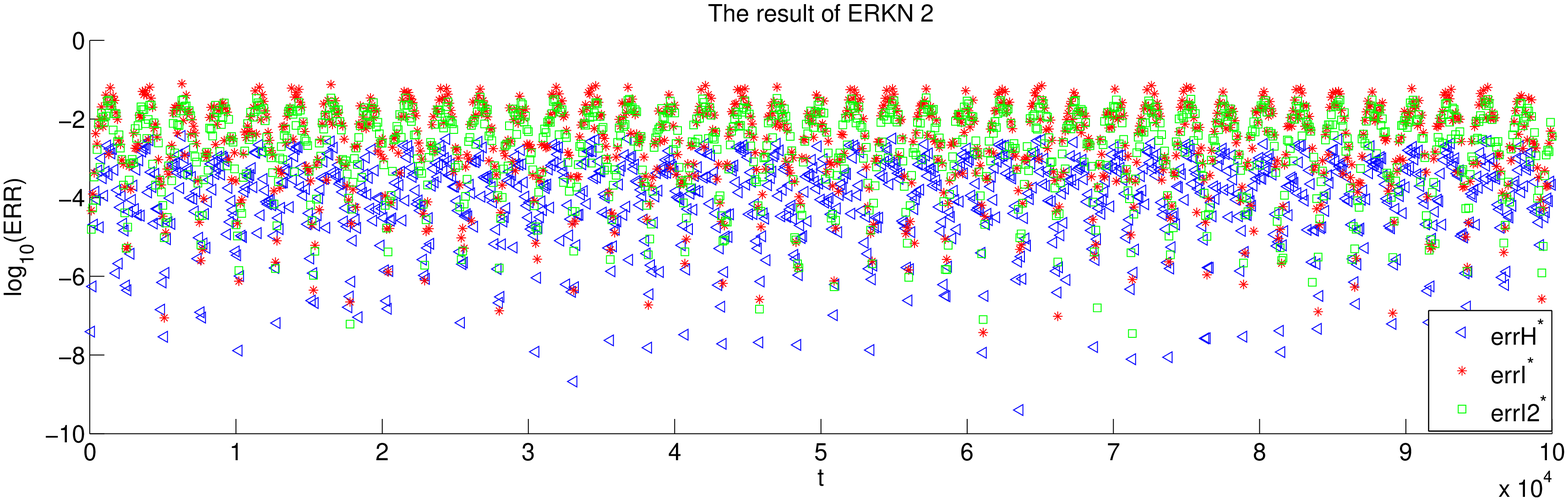}
\end{tabular}
\caption{The errors  of energies (up) and modified energies (down) against $t$   for ERKN 2.}%
\label{fig2}%
\end{figure}

\begin{figure}[ptb]
\centering\tabcolsep=2mm
\begin{tabular}
[c]{ccc}%
\includegraphics[width=13cm,height=3.5cm]{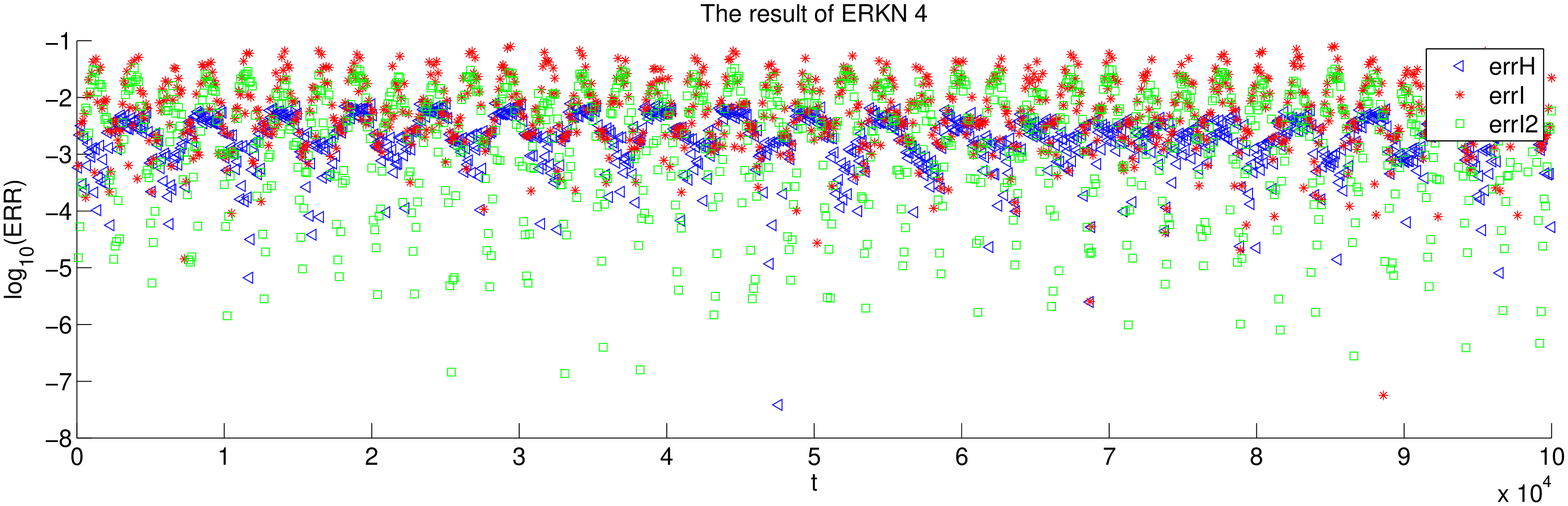}  \\
\includegraphics[width=13cm,height=3.5cm]{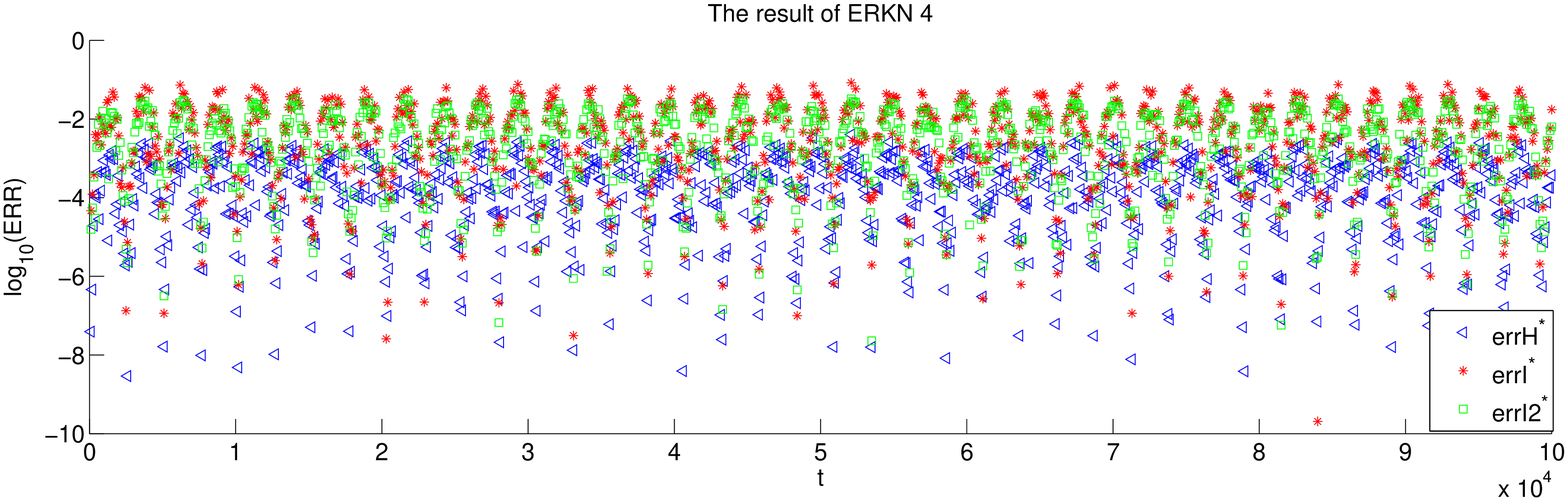}
\end{tabular}
\caption{The errors  of energies (up) and modified energies (down)  against $t$  for ERKN 4.}%
\label{fig3}%
\end{figure}

The system is integrated in the interval $[0,10000]$ with $h=0.01$.
First the errors of the total energy $H$ and oscillatory energy $I$
and $I_2$ against $t$ for  ERKN3   and the errors of the modified
energies $H^{*},I^{*},I^{*}_2$ for ERKN1 are shown in Fig.
\ref{fig1}. Then we present the   energies
 and the   modified energies for ERKN 2 and 4 in Figs.
\ref{fig2} and \ref{fig3}, respectively.

 From the numerical results, it follows that the non-symmetric ERKN1 can not preserve
the energies and  the symmetric and symplectic ERKN3 can
approximately  conserve the energies very well over long times. For
the symmetric ERKN 2 and   4 not satisfying the  condition
\eqref{new-cond}, they approximately conserve the modified energies
better than the original energies.

\section{Conclusions} \label{sec:conclusions}

This paper studied  the long-time behaviour    of   ERKN integrators
 for muti-frequency
highly oscillatory Hamiltonian systems. The modulated
multi-frequency Fourier expansion of ERKN integrators were developed
and by which,   we showed the long-time numerical energy
conservation of the integrators. Our next work  will be devoted to
the long-time analysis of ERKN integrators for multi-frequency
highly oscillatory Hamiltonian systems under minimal non-resonance
conditions.

\section*{Acknowledgements}

The first author is grateful to Professor Christian Lubich for the
helpful  discussions on the topic of modulated Fourier expansions.

The research is supported in part  by the Alexander von Humboldt
Foundation, by the  Natural Science Foundation of Shandong Province
(Outstanding Youth Foundation) under Grant ZR2017JL003, and by the
National Natural Science Foundation of China under Grant 11671200.

\end{document}